\documentclass[12pt]{amsart}
\usepackage[margin=2.9cm]{geometry}
\usepackage[foot]{amsaddr}
\usepackage{amsmath}
\usepackage[shortlabels]{enumitem}
\usepackage[center]{caption}
\usepackage{latexsym}
\usepackage{amsfonts}
\usepackage{graphicx}
\usepackage{caption}
\usepackage{todonotes}
\usepackage{amssymb, amsthm,latexsym, mathdots, mathtools, graphicx, mathtools, xcolor, xkeyval, tikz, tkz-euclide, adjustbox, float, array, subcaption, listings, comment}
\usepackage[colorlinks,linkcolor=blue!60!black,citecolor=black]{hyperref}
\usepackage{cleveref}
\renewcommand{\qedsymbol}{$\blacksquare$}
\newtheorem{theorem}{Theorem}[section]
\newtheorem{lemma}[theorem]{Lemma}

\newtheorem{observation}[theorem]{Observation}
\newtheorem{corollary}[theorem]{Corollary}
\newtheorem{claim}{Claim}
\newtheorem{case}{Case}

\newcommand{\claw}{\mbox{claw}}
\newcommand{\paw}{\mbox{paw}}
\newcommand{\diam}{\mbox{diamond}}
\newcommand{\trig}{\mbox{triangle}}

\newcounter{tbox}
\newcommand{\clm}[1]{\vspace{0.1cm}\medskip\refstepcounter{tbox}\noindent{\parbox{\textwidth}{(\thetbox) \emph{#1}}}\vspace{0.1cm}}

\title[Non-empty intersection of longest paths]{Non-empty intersection of longest paths\\ in $P_5$-free and claw-free graphs}
\author{Paloma T. Lima$^{\mathsection}$}
\author{Amir Nikabadi$^{\mathsection}$}
\address{$\mathsection$ IT University of Copenhagen, Denmark}
\address{Lima and Nikabadi acknowledge the support of the Independent Research Fund Denmark grant agreement number 2098-00012B.}
\email{\{palt,amir\}@itu.dk}
\date{\today}

\begin{document}
\maketitle

\begin{abstract}
A family $\mathcal{F}$ of graphs is a \textit{Gallai family} if for every connected graph $G\in \mathcal{F}$, all longest paths in $G$ have a common vertex. While it is not known whether $P_5$-free graphs are a Gallai family, Long Jr., Milans, and Munaro [The Electronic Journal of Combinatorics, 2023] showed that this is \emph{not} the case for the class of claw-free graphs. We give a complete characterization of the graphs $H$ of size at most five for which $(\text{claw}, H)$-free graphs form a Gallai family. We also show that $(P_5, H)$-free graphs form a Gallai family if $H$ is a triangle, a paw, or a diamond. Both of our results are constructive.
\end{abstract}

\section{Introduction}

In 1966, Gallai~\cite{Gallai68} proposed the following question: \textit{do all longest paths in a connected graph share at least one vertex?} The question was answered in the negative by a counterexample due to Walther~\cite{walther1969nichtexistenz} (a minimal counterexample was later provided by Zamfirescu~\cite{zamfirescu1976longest}).
On the other hand, multiple graph classes were shown to have a positive answer to Gallai's question. 
We say a vertex in a graph is a \textit{Gallai vertex} if it belongs to all of its longest paths.
We say a graph class $\mathcal{F}$ forms a \emph{Gallai family} if for every connected graph $G\in \mathcal{F}$, $G$ has a Gallai vertex.
Examples of classes of graphs that have been shown to be Gallai families include circular-arc graphs~\cite{BalisterGLS04,Joos15}, dually chordal graphs~\cite{JobsonKLW16}, bipartite permutation graphs~\cite{CERIOLI2020111717}, graphs of treewidth at most 2~\cite{ChenEFHSYY17}, $P_4$-sparse graphs~\cite{cerioli2020intersection} and graphs of matching number at most 3~\cite{Chen15}. It remains a challenging question to understand which structural properties of a graph ensure the intersection of all its longest paths is non-empty.

Graph classes defined by a finite family of forbidden induced subgraphs have also been investigated. In particular, $2P_2$-free graphs~\cite{golanshan2k2} and $(P_5,K_{1,3})$-free graphs~\cite{cerioli2020intersection} were shown to form Gallai families. More recently, Long Jr., Milans, and Munaro \cite{long2023non} initiated a systematic study of the problem on $H$-free graphs. They proved that, if $H$ is a graph on at most four vertices, connected $H$-free graphs form a Gallai family if and only if $H$ is a linear forest. In particular, this implies that claw-free graphs do \emph{not} form a Gallai family. Gao and Shan~\cite{gao2021nonempty} identified multiple subclasses of claw-free graphs that have a positive answer to Gallai's question, including $(claw,P_6)$-free graphs.  

In this work, motivated by these recent results, we continue the systematic investigation of graph classes defined by two forbidden induced subgraphs, with a focus on subclasses of claw-free graphs and $P_5$-free graphs. We first prove the following theorem, which together with the results of Gao and Shan~\cite{gao2021nonempty} and Long Jr., Milans, and Munaro \cite{long2023non}, implies a complete classification of $(claw,H)$-free graphs that form a Gallai family, for $H$ of size at most five.

\begin{theorem}\label{thm:main-claw,Hinintgro}
Let $H$ be one of the graphs $P_3 + 2P_1, K_3+2P_1, 2P_2+P_1,$ or $P_2+3P_1$. Let $\mathcal{C}$ be the class of graphs restricted to $(claw, H)$-free graphs, then $\mathcal{C}$ forms a Gallai family.
Moreover, there exists a polynomial-time algorithm which finds Gallai vertices on graphs in $\mathcal{C}$. 
\end{theorem}

We then shift our attention to $P_5$-free graphs. It remains an open problem whether $P_5$-free graphs form a Gallai family. We conjecture this to be the case. We provide evidence in this direction with the following result.

\begin{theorem}\label{thm:main}
Let $H$ be one of the graphs triangle, paw, or diamond.
	Let $\mathcal{G}$ be the class of graphs restricted to $(P_5, H)$-free graphs, then $\mathcal{G}$ forms a Gallai family. Moreover, there exists a polynomial-time algorithm which finds Gallai vertices on graphs in $\mathcal{G}$. 
\end{theorem}

Theorem~\ref{thm:main}, together with the results from Long Jr., Milans, and Munaro \cite{long2023non} leaves open only two cases ($H=K_4$ and $H=C_4$) towards showing that for any graph $H$ of size at most four, $(P_5,H)$-free graphs form a Gallai family. 

\subsection*{A Short Background.} While here we focus on classifying whether a graph class is a Gallai family, other works have sought to determine, for a graph $G$, what is the size of the smallest set $S\subset V(G)$ such that every longest path of $G$ intersects $S$. Such a set is called a \emph{longest path transversal}. Note that a graph class $\mathcal{F}$ is a Gallai family if all the graphs in $\mathcal{F}$ admit longest path transversals of size one. In 2021, Long Jr., Milans and Munaro~\cite{sublinearLPT} showed that any graph admits a longest path transversal of sublinear size, improving on an earlier bound by Rautenbach and Sereni~\cite{rautenbach2014}. Some graph classes have been shown to admit longest path transversals of constant size (with a constant larger than one). This is the case for $(P_5,cricket)$-free graphs~\cite{cerioli2020intersection} and graphs of constant treewidth~\cite{rautenbach2014,CERIOLI2020111717}. 

\section{Preliminaries}\label{sec:prelim}
Graphs in this paper are finite and simple. Let $G = (V(G),E(G))$ be a graph. A \textit{clique} in $G$ is a set of pairwise adjacent
vertices. A \textit{stable set} or \textit{independent set} in $G$ is a set of pairwise non-adjacent vertices. A \textit{linear forest} is a forest whose components are paths. We let $P:= p_1p_2\dots p_k$ denote a path in $G$. We call the vertices $p_1$ and $p_k$ the endpoints of $P$ and say that $P$ is a path from $p_1$ to $p_k$. For a vertex $x\in V(G)$, we say $x$ is \textit{between $p_i$ and $p_j$}, if $x$ is in the path $p_i \dots p_j$. We also say \textit{$x$ has no neighbor from $p_i$ to $p_j$}, if there is no vertex $y$ between $p_i$ and $p_j$ such that $xy$ is an edge.
A path (cycle) in a graph is \textit{Hamiltonian}, if it contains all vertices of the graph.
A graph is \textit{$k$-connected} if it has at least $k+1$ vertices and no vertex separator of size at most $k-1$. A $diamond$ is a graph obtained from $K_4$ by removing an edge. A \textit{triangle} is a $K_3$.
If $X \subseteq V(G)$ we denote the subgraph induced on $X$ by $G[X]$. For disjoint $X, Y \subseteq G$, we say that $X$ is \textit{complete} to $Y$ if every vertex in $X$ is adjacent to every vertex in $Y$, and $X$ is
\textit{anticomplete} to $Y$ if there are no edges between $X$ and $Y$.
For given graphs $G$ and $H$, we say that $G$ is \textit{$H$-free} if $G$ does not contain $H$ as an induced subgraph. We say $G$ is $(H_{1}, H_{2})$-free if $G$ do not contain $H_1$ and $H_2$ as an induced subgraphs. We let $P_n$, $C_n$, and $K_n$ denote the chordless path, chordless cycle, and the complete graph on $n$ vertices. For integer $t \geq 1$, we denote by $tP_n$ the graph obtained from the disjoint union of $t$ copies of the $n$-vertex path, and for graphs $G_1$, $G_2$, we write $G_1 + G_2$ to denote the disjoint union of $G_1$ and $G_2$. A subset $D \subseteq V(G)$ is a \textit{dominating set} if each vertex of $G$ either belongs to $D$ or is adjacent to some vertex of $D$.

\begin{observation}[\cite{cerioli2020intersection}]\label{obs:domupperbound}
    Let $D\subseteq V(G)$ be a dominating set and let $\mathcal{S} \subseteq V(G)$ be the minimum size of a longest path transversal in $G$. Then, $|\mathcal{S}|\leq |D|$.
\end{observation}

\begin{observation}\label{obs:ifp1pkinE}
    Let $G$ be a connected graph and $P:=p_1\ldots p_k$ be a longest path of $G$. If $p_1p_k\in E(G)$, then $P$ is a Hamiltonian path.
\end{observation}

\begin{proof}[Proof of \Cref{obs:ifp1pkinE}]
Suppose $P$ is not a Hamiltonian path. Since $G$ is connected, there exists $x\notin V(P)$ that has a neighbor in $V(P)$. Let $p_i$ be such a neighbor. Then $xp_i\ldots p_kp_1\ldots p_{i-1}$ is a path longer than $P$ in $G$ which is a contradiction.
\end{proof}

The following results (see also \cite{golanshan2k2} for the case of $2P_2$-free graphs) will partly play a role in our results:

\begin{theorem}[\cite{long2023non}]\label{necessityforlinearforests}
Let $H$ be a graph. If the class of connected $H$-free graphs forms a Gallai family, then $H$ is a linear forest on at most 9 vertices.
\end{theorem}

\begin{theorem}[\cite{long2023non}]\label{2p2-free}
Let $H$ be a graph on at most four vertices. The class of connected $H$-free graphs forms a Gallai family if and only if $H$ is a linear forest.
\end{theorem}

\section{$(\claw,H)$-free graphs}\label{sec:claw,H}

We begin the section with two observations. 
The graph $\mathcal{B}$ depicted in Figure~\ref{fig:scatt} is the example due to Walther~\cite{walther1969nichtexistenz} answering Gallai's question in the negative. In particular, the graph $\mathcal{B}$ is constructed from the Petersen graph (which has no Hamiltonian cycle) by blowing up an arbitrary vertex into a set of three vertices, each of degree 1.
 
Long Jr., Milans and Munaro~\cite{long2023non} showed that one can construct a claw-free graph $\mathcal{B^{+}}$ from the graph $\mathcal{B}$ (by replacing each vertex of degree three in $\mathcal{B}$ with a triangle) such that $\mathcal{B^{+}}$ has no Gallai vertex. This implies that:

\begin{observation}\label{obs:counterexample}
The class of claw-free graphs does not form a Gallai family.
\end{observation}

Second, Gao and Shan \cite{gao2021nonempty} proved that the family of $(H_1,H_2)$-free graphs such that $H_1$ and $H_2$ are connected and every 2-connected $(H_1,H_2)$-free graph has a Hamiltonian cycle, forms a Gallai family. Let us state this precisely. For integers $i, j, k\geq 0$, we denote by $N_{i,j,k}$ the graph obtained from a triangle by appending disjoint paths of length $i, j, k$ at each vertex of the triangle. A \textit{paw} is a $N_{1,0,0}$, and a \textit{bull} is a $N_{1,1,0}$.

\begin{theorem}[\cite{gao2021nonempty}]\label{thm:dich-for-claw-H}
Let $0\leq i \leq 3$. The class of $(claw,H)$-free graphs forms a Gallai family if $H$ is one
of the graphs $P_4, P_5, P_6, bull$, $N_{2,1,0}$, or $N_{i,0,0}$.
\end{theorem}

Observe that the classes considered in Theorem~\ref{thm:main-claw,Hinintgro} are not contained in any of the classes covered by Theorem~\ref{thm:dich-for-claw-H}. 
The rest of the section is then devoted to proving Theorem~\ref{thm:main-claw,Hinintgro}. We will need the following observation for our proofs in this section:

\begin{figure}[t]
    \centering
    \includegraphics[scale=1]{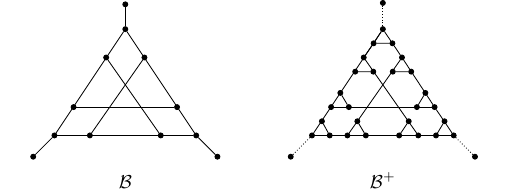}
    \caption{Left: The graph $\mathcal{B}$ with no Gallai vertex. Right: A claw-free graph $\mathcal{B}^{+}$ obtained from $\mathcal{B}$ (dotted lines show paths of length at least one).}
   \label{fig:scatt}
\end{figure}

\begin{observation}\label{obs:bridge1}
Let $G$ be a claw-free graph, $P=p_1p_2\ldots p_k$ be a longest path of $G$, and let $u\in V(G)$ be such that $u\notin V(P)$. Then the following holds.
\begin{itemize}
\item If $up_i\in E(G)$, then $p_{i-1}p_{i+1}\in E(G)$;
\item Moreover, if $u$ has a neighbor $v\notin V(P)$, then $\{u,v,p_{i-1},p_{i+1}\}$ induces a $2P_2$ in $G$.
\end{itemize}
\end{observation}

\begin{proof}[Proof of \Cref{obs:bridge1}]
For the first part of Observation~\ref{obs:bridge1}; see that $up_{i-1}\notin E(G)$, otherwise $p_1\ldots p_{i-1}up_i\ldots p_k$ is a path longer than $P$. The same holds for the edge $up_{i+1}$. Hence we conclude $p_{i-1}p_{i+1}\in E(G)$, otherwise $\{p_{i-1}, p_i, p_{i+1}, u\}$ induces a claw in $G$, a contradiction. For the second part; suppose $u$ has a neighbor $v\notin V(P)$. If $p_{i-1}v\in E(G)$, then $p_1\ldots p_{i-1}vup_{i}\ldots p_k$ is a path longer than $P$, again, a contradiction. The same holds if $p_{i+1}v\in E(G)$.
\end{proof}
\begin{theorem}\label{thm:claw,p3+2p1}
The class of $(\claw,P_3+2P_1)$-free graphs forms a Gallai family.
\end{theorem}

\begin{proof}
Let $G$ be a connected $(\claw,P_3+2P_1)$-free graph. If $G$ does not have an induced subgraph isomorphic to $P_3+P_1$, by Theorem~\ref{2p2-free}, $G$ has a Gallai vertex. Hence, let $D\subset V(G)$ be a set inducing a $P_3+P_1$. Note that $D$ must be a dominating set, otherwise $G$ would have an induced $P_3+2P_1$. Let $D=\{x,y,z,w\}$, where $xyz$ is a $P_3$ and $w$ is an isolated vertex in~$D$. Let $P:= p_1p_2 \ldots p_k$ be a longest path in $G$.

\begin{claim}\label{clm:xyz-are-in-P}
$V(P) \cap \{x,y,z\} \neq \emptyset$.
\end{claim}

\begin{proof}[Proof of \Cref{clm:xyz-are-in-P}]\renewcommand{\qedsymbol}{}
Suppose not. First, observe that $p_1$ and $p_k$ are not adjacent to any vertex of $\{x,y,z\}$, otherwise one can extend $P$ into a longer path. Since $D$ is a dominating set in $G$, $p_1$ and $p_k$ are both adjacent to $w$. Note that this implies $w \in V(P)$, otherwise we can extend $P$ to $w$. Consider the set $\{x,y,z, p_1, p_k\}$. Since $G$ is ($P_3+2P_1$)-free,  $p_1p_k \in E(G)$. Furthermore, as $G$ is connected,
there must exist a (shortest) path between $\{x,y,z\}$ and $V(P)$. Let $p_i$ be the vertex in $V(P)$ that is closest to $\{x,y,z\}$ and let $Q$ be a shortest path between $\{x,y,z\}$ and $p_i$. Let $q$ be the vertex in $Q$ that is adjacent to $p_i$. Note that it could be that $q\in\{x,y,z\}$. Then $ q p_i \dots p_kp_1 \dots p_{i-1}$ is a longer path, a contradiction. This proves \Cref{clm:xyz-are-in-P}.
\end{proof}

We proceed by showing that $y$ is a Gallai vertex. Assume for the sake of contradiction that $y \notin V(P)$, then by \Cref{clm:xyz-are-in-P} either of $x$ or $z$ must be in $V(P)$. By symmetry we may assume $x\in V(P)$. Let $a$ and $b$ be the vertices appearing before and after $x$ in $P$, respectively. Since $y\notin V(P)$, by Observation~\ref{obs:bridge1}, $ab \in E(G)$. 
Clearly, $p_1y$ cannot be an edge, otherwise we can extend $P$ to $y$. Moreover, $p_1x\notin E(G)$, otherwise $yxp_1\dots ab \dots p_k$ is a longer path in $G$. In a similar way, we conclude $p_1z\notin E(G)$. Indeed, if that was the case, $z$ would be in $P$, and the same arguments would yield a longer path. The same holds for $p_k$: it cannot be adjacent to any vertex in $\{x,y,z\}$. Consider the set $\{x,y,z,p_1,p_k\}$. Since $G$ is $(P_3+2P_1)$-free, we have $p_1p_k\in E(G)$. Then $yxb\ldots p_kp_1\ldots a$ is a path longer than $P$ in $G$, a contradiction.
Hence, $y$ is common to all longest paths in $G$. This finishes the proof of Theorem~\ref{thm:claw,p3+2p1}.
\end{proof}

\begin{observation}\label{alg:claw,P3+2P1}
    Let $G$ be a connected $(\claw, P_3+2P_1)$-free graph. There exists a polynomial-time algorithm which finds a Gallai vertex in $G$.
\end{observation}

\begin{proof}[Proof of Observation~\ref{alg:claw,P3+2P1}]{}
Consider the following simple algorithm: in polynomial time we can check whether $G$ contains an induced copy of $P_3+P_1$. If the answer is yes, that is, if $G$ contains $P_3+P_1$ as an induced subgraph, then we return the middle vertex of the $P_3$, say $v$. It follows from the proof of Theorem~\ref{thm:claw,p3+2p1} that $v$ is a Gallai vertex. If the answer is no, that is, if $G$ is $(\claw, P_3+P_1)$-free, then by the result of \cite{long2023non} (see Theorem 6 in \cite{long2023non}), every vertex of degree at least $\Delta(G)-1$ is a Gallai vertex, and such vertex can be found in polynomial time.
\end{proof}

\medskip

A similar proof as above works for the case in which $G$ is a $(\claw, K_3+2P_1)$-free graph, in which case we would have that $\{x,y,z\}$ induces a triangle instead. Hence, we have the following:

\begin{theorem}\label{thm:claw,K3+2P1}
The class of $(\claw, K_3+2P_1)$-free graphs forms a Gallai family. 
\end{theorem}

To apply the argument of Observation~\ref{alg:claw,P3+2P1} for Theorem~\ref{thm:claw,K3+2P1}, note that the class of $(\claw, N_{3,0,0})$-free graphs contains the class of $(\claw, K_3+P_1)$-free graphs. The former class is a Gallai family, by \cite{gao2021nonempty}. We leave it to the reader to check that the proof of \cite{gao2021nonempty} for $(\claw, N_{3,0,0})$-free graphs is constructive.

\medskip

Next, we prove the following:

\begin{theorem}\label{thm:claw,2P_2+P_1}
    The class of $(\claw, 2P_2+P_1)$-free graphs forms a Gallai family.
\end{theorem}

\begin{proof}
Let $G$ be a connected $(\claw,2P_2+P_1)$-free graph. If $G$ does not have an induced~$2P_2$, then it follows by Theorem~\ref{2p2-free} that $G$ has a Gallai vertex. Hence, let $\{x,y,z,w\}\subset V(G)$ be a set inducing a $2P_2$, where $xy\in E(G)$ and $zw\in E(G)$. We may assume $d(x)\geq d(y)$. We will now show that $x$ is a Gallai vertex in $G$. Let $P=p_1\ldots p_k$ be a longest path in $G$ and suppose for a contradiction that $x \notin V(P)$. We first claim the following.

\begin{claim}\label{claim:edgeintP}
$V(P)\cap\{x,y\}\neq\emptyset$.
\end{claim}

\begin{proof}[Proof of \Cref{claim:edgeintP}]\renewcommand{\qedsymbol}{}
Suppose this is not the case. Let $y'\in V(G)$ be the vertex in $V(P)$ that is closest to $\{x,y\}$ and let $Q$ be a shortest path between $\{x,y\}$ and $y'$. Let $y''$ be the vertex adjacent to $y'$ in $Q$. Note that it may be that $y''\in\{x,y\}$.
Let $a$ and $b$ be the vertices before and after $y'$ in $P$, respectively. Consider the set $\{a,b,y',y''\}$. By Observation~\ref{obs:bridge1}, $ab\in E(G)$ and note that $\{ab, xy\}$ induces a $2P_2$ in $G$. Indeed, if $ay\in E(G)$ then $p_1\ldots a y\cdot Q\cdot y'b\ldots p_k$ is a path longer than $P$. A similar argument shows that $ax,by,bx\notin E(G)$. Consider the set $\{a,b,x,y,p_1\}$. Since $G$ is $(2P_2+P_1)$-free, $p_1$ must be adjacent to a vertex in $\{a,b,x,y\}$. Since, $x$ and $y$ are not in $P$, neither of $p_1x$ nor $p_1y$ is an edge, otherwise we could extend $P$ to $x$ or $y$. If $p_1b\in E(G)$, then $p_k\ldots bp_1\ldots y'y''$ is a path longer than $P$ in $G$. Hence, we have $p_1a\in E(G)$. Note that $\{p_1a, xy\}$ induces a $2P_2$. Consider the set $\{p_1, a, x, y, p_k\}$. Now, $p_k$ must be adjacent to a vertex in $\{p_1,a,x,y\}$. Similarly to $p_1$, $p_k$ cannot be adjacent to $x$ and $y$. If $p_kp_1\in E(G)$, then $y''y'\ldots p_kp_1\ldots a$ is a path longer than $P$. Hence $p_ka\in E(G)$. Let $a'$ be the vertex before $a$ in $P$. Then $a'\ldots p_1ap_k\ldots y'y''$ is a path longer than $P$, a contradiction. This concludes the proof of \Cref{claim:edgeintP}.
\end{proof}

Observe that, since $x\notin V(P)$, by \Cref{claim:edgeintP}, we have $y\in V(P)$. Moreover, if $a$ and $b$ are the vertices before and after $y$ in $P$, by Observation~\ref{obs:bridge1}, we have $ab\in E(G)$. 

\begin{claim}\label{claim:2k2in2k2+p1}
$\{x,y,p_1,p_2\}$ induces a $2P_2$ in $G$.
\end{claim}

\begin{proof}[Proof of \Cref{claim:2k2in2k2+p1}]\renewcommand{\qedsymbol}{}
Note that $xp_1\notin E(G)$, otherwise $P$ could be extended to $x$. If $xp_2\in E(G)$, by Observation~\ref{obs:bridge1}, $p_1p_3\in E(G)$. Hence, $xp_2p_1p_3\ldots p_k$ is a path longer than $P$. If $yp_1\in E(G)$, then $xyp_1\ldots ab\ldots p_k$ is a path longer than $P$. Finally, assume $yp_2\in E(G)$. Since $d(x)\geq d(y)$, and $\{a,b\}\subset N(y)\setminus N(x)$, $x$ must have a neighbor $x'$ that is not a neighbor of $y$. If $x'\notin V(P)$, then $x'xyp_2\ldots ab\ldots p_k$ is a path longer than $P$, since both $x$ and $x'$ were added, while only $p_1$ was excluded. Thus, $x'\in V(P)$. We may assume, by symmetry, that $x'$ lies between $b$ and $p_k$ in $P$, as in this claim, the same argument would work for the edge $p_{k-1}p_k$ instead of $p_1p_2$. Let $c$ and $d$ be the vertices before and after $x'$ in $P$, respectively. By Observation~\ref{obs:bridge1}, $cd\in E(G)$ and $\{x,y,c,d\}$ induces a $2P_2$ in $G$. Consider the set $\{x,y,c,d,p_1\}$. We already argued $p_1x,p_1y\notin E(G)$. If $p_1d\in E(G)$, then $p_k\ldots dp_1\ldots x'x$ is a path longer than $P$. If $p_1c\in E(G)$, then $p_3\ldots ab\ldots cp_1p_2yxx'\ldots p_k$ is a path longer than $P$, a contradiction. Hence we conclude that $yp_2\notin E(G)$ and therefore \Cref{claim:2k2in2k2+p1} holds.
\end{proof}

Now consider the set $\{x,y,p_1,p_2,p_k\}$. By \Cref{claim:2k2in2k2+p1}, $\{x,y,p_1,p_2\}$ induces a $2P_2$. If $xp_k\in E(G)$, $P$ could be extended to $x$. If $yp_k\in E(G)$, then $xyp_k\ldots ba\ldots p_1$ is a path longer than $P$. By Observation~\ref{obs:ifp1pkinE}, $p_kp_1\notin E(G)$. Then we must have $p_kp_2\in E(G)$, otherwise $\{x,y,p_1,p_2,p_k\}$ induces a $2P_2+P_1$ in $G$. Now we use again an argument similar to the one in the proof of \Cref{claim:2k2in2k2+p1}, where we had the edge $yp_2$. This time, instead we will use the fact that $p_kp_2\in E(G)$.
Since $d(x)\geq d(y)$, and $\{a,b\}\subset N(y)\setminus N(x)$, $x$ must have a neighbor $x'$ that is not a neighbor of $y$. If $x'\notin V(P)$, then $x'xyb\ldots p_kp_2\ldots a$ is a path longer than~$P$. Hence, $x'\in V(P)$. We may assume $x'$ lies between $p_1$ and $a$ in $P$ (if not, a similar argument would work with the edge $p_{k-1}p_k$).
Let $c$ and $d$ be the vertices before and after $x'$ in $P$, respectively. Consider the set $\{x,y,c,d,p_1\}$. By Observation~\ref{obs:bridge1}, $cd\in E(G)$. Note that $\{x,y,c,d\}$ induces a $2P_2$ in $G$. Similarly, as before, $p_1x,p_1y\notin E(G)$, otherwise longer paths could be obtained. If $p_1c\in E(G)$, then $p_3\ldots cp_1p_2p_k\ldots x'x$ is a path longer than $P$. So $p_1c\notin E(G)$. Then we must have $p_1d\in E(G)$, otherwise $\{x,y,c,d,p_1\}$ induces a $2P_2+P_1$. To conclude the proof, consider the set $\{x,y,c,d,p_k\}$. Similarly to $p_1$, we cannot have edges between $p_k$ and $x$ or $y$. Let $c'$ be the vertex before $c$ in $P$ and $d'$ be the vertex after $d$. If $p_kc\in E(G)$, then $xx'cp_k\ldots dp_1\ldots c'$ is a path longer than $P$ in $G$. If $p_kd\in E(G)$, then $d'\ldots p_kdp_1\ldots x'x$ is a path longer than $P$. Therefore we conclude $\{x,y,c,d,p_k\}$ induces a $2P_2+P_1$ in $G$, a contradiction. This implies that $x$ must belong to any longest path $P$ of $G$, hence finishes the proof of Theorem~\ref{thm:claw,2P_2+P_1}.
\end{proof}

\medskip

\begin{observation}\label{alg:claw,2P2+P1}
      Let $G$ be a connected $(\claw, 2P_2+P_1)$-free graph. There exists a polynomial-time algorithm which finds a Gallai vertex in $G$.
\end{observation}

The proof of Observation~\ref{alg:claw,2P2+P1} is is quite similar to the proof of Observation~\ref{alg:claw,P3+2P1}. We remind the reader that in a connected $2P_2$-free graph, every vertex of maximum degree is common to all longest paths (see Theorem 1 in \cite{golanshan2k2}).

\medskip

The remainder of this section is devoted to settling Theorem~\ref{thm:main-claw,Hinintgro} for $(\claw, P_2+3P_1)$-free graphs. This will complete our picture for Gallai families on $(\claw, H)$-free graphs with $|H|\leq 5$. Although the big picture behind the proof is
similar to the proof of Theorem~\ref{thm:claw,p3+2p1} to Theorem~\ref{thm:claw,2P_2+P_1}, we are compelled to take extra steps to mold the details. Recall that a \textit{cut vertex} $x \in V(G)$ of $G$ is a vertex such that $G-x$ is disconnected. 

\begin{figure}[t]
    \centering
    \includegraphics[scale=1.5]{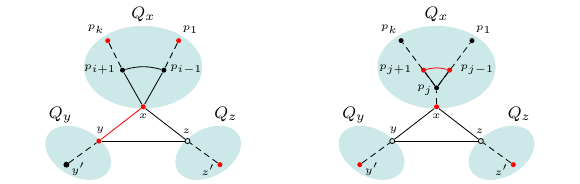}
    \caption{Proof of Lemma~\ref{lem:three-cut-vertex}. Left: the case $x\in V(P)$. Right: the case $x\notin V(P)$. Dashed lines show paths of
arbitrary length.}
   \label{fig:proofoflemma}
\end{figure}

\begin{lemma}\label{lem:three-cut-vertex}
    Let $G$ be a connected $(\claw, P_2+3P_1)$-free graph with three cut vertices $\{x,y,z\} \subseteq V(G)$ inducing a triangle in $G$. Then $x$, $y$ and $z$ are Gallai vertices.
\end{lemma}

\begin{proof}
For $\{x,y,z\} \subseteq V(G)$ we denote by $Q_{x}$ a subset of vertices of $G$ that are reachable from $x$ by a path, say $\mathcal{L}$, such that $y,z\notin V(\mathcal{L})$. $Q_{y}$ and $Q_z$ are defined similarly. Observe that since $\{x,y,z\}$ are cut vertices inducing a triangle; $xy$, $xz$, and $yz$ are the only edges between $Q_{x}, Q_{y},Q_{z}$. Let $P:=p_1\dots p_k$ be a longest path in $G$. It is straightforward to check that if $P$ contains two vertices of $\{x,y,z\}$, it must also include the third one and so the claim holds. Let $q \in \{x,y,z\}$ and suppose (without loss of generality) that $q = x$. Assume for a contradiction that $V(P) \subseteq G[Q_x]$. There are two possibilities:

\begin{itemize}\setlength\itemsep{0.5em}
    \item $x\in V(P)$: let $x = p_i$, $1< i < k$. Let $z'$ be a vertex in $G[Q_z]$ such that $z'\neq z$. See Figure~\ref{fig:proofoflemma}, left. We claim that $\{x, y, p_1, p_k, z'\}$ induces a $P_2+3P_1$ in $G$. Since $y,z \in N(x)\setminus V(P)$, by Observation~\ref{obs:bridge1} we have $p_{i-1}p_{i+1} \in E(G)$. Then $p_1x \notin E(G)$, otherwise $z' \dots zyxp_1 \dots p_{i-1}p_{i+1} \dots p_k$ is a longer path, symmetrically, $p_kx \notin E(G)$. Moreover, since $x$ is a cut vertex and $y,z\notin G[Q_x]$, neither of $yp_1$ nor $yp_k$ is an edge. Note that since $z$ is a cut vertex $\{x,y,z'\}$ induces a $P_2+P_1$, also $p_1z', p_{k}z' \notin E(G)$. This combined with the fact that $p_1p_k \notin E(G)$, which follows from Observation~\ref{obs:ifp1pkinE}, shows $\{x, y, p_1, p_k, z'\}$ indeed induces a $P_2+3P_1$ in $G$, a contradiction.

    \item $x \notin V(P)$: let $p_j$, $1<j<k$, be the vertex in $V(P)$ that is closest to $x$ and let $\hat{P}$ be a shortest path between $x$ and $p_j$. See Figure~\ref{fig:proofoflemma} to the right. Let $x'$ be the vertex adjacent to $p_j$ in $\hat{P}$. Note that it may be that $x'=x$. Again, observe that $j\neq 1$ and $j\neq k$, otherwise we could extend $P$ to $xyz$. Since $x \notin V(P)$ we have (i) by Observation~\ref{obs:bridge1} that $p_{j-1}p_{j+1}\in E(G)$, and (ii) neither of $xp_{j-1}$ nor $xp_{j+1}$ is an edge. But this along with the fact that $x$ is a cut vertex, implies that $\{p_{j-1}, p_{j+1}, x, y', z'\}$ induces a $P_2+3P_1$, where $y'$ be a vertex in $G[Q_y]$ such that $y'\neq y$ and $z'$ is defined analogously. This is a contradiction.  
\end{itemize}

This completes the proof of Lemma~\ref{lem:three-cut-vertex}.
\end{proof}

\begin{theorem}\label{thm:claw,p2+3p1}
    The class of $(\claw, P_2+3P_1)$-free graphs forms a Gallai family.
\end{theorem}

\begin{proof}
Let $G$ be a connected $(\claw, P_2+3P_1)$-free graph. By Lemma~\ref{lem:three-cut-vertex}, we may assume that $G$ has no three cut vertices inducing a triangle. If $G$ is a $(P_2+2P_1)$-free graph, then by Theorem~\ref{2p2-free}, $G$ has a Gallai vertex. Hence, we may assume $G$ has a set $D\subset V(G)$ inducing a $P_2+2P_1$. Let $D = \{x,y, z,w\}$, with $xy\in E(G)$. Note that since $G$ is $(P_2+3P_1)$-free, $D$ must be a dominating set in $G$. We may assume, without loss of generality, that $d_{G}(x)\geq d_{G}(y)$. We wish to show that $x$ is a Gallai vertex in $G$. As before, let $P:= p_1 \dots p_k$ be a longest path of $G$ and suppose for the sake of contradiction that $x\notin V(P)$.
We begin with a sequence of claims.

\begin{claim}\label{clm:intersection-geq-2}
    $|V(P) \cap V(D)| \geq 2$
\end{claim}

\begin{proof}[Proof of \Cref{clm:intersection-geq-2}]\renewcommand{\qedsymbol}{}
Suppose not. Since $D$ is a dominating set, by Observation~\ref{obs:domupperbound}, $V(P) \cap V(D) \neq \emptyset$. Assume first that $V(P) \cap V(D) = \{y\}$. Since $D$ is a dominating set, we have $p_1y \in E(G)$, and $p_ky \in E(G)$, otherwise it would be possible to extend $P$. Let $y'$ and $y''$ be the vertices adjacent to $y$ in $P$. Observe that as $x \notin V(P)$, the edge $y'y'' \in E(G)$ by Observation~\ref{obs:bridge1}. But then $xyp_1\dots y'y'' \dots p_k$ is a longer path, a contradiction. The case $V(P) \cap V(D) = \{x\}$ is symmetrical. Next, assume $V(P) \cap V(D) = \{z\}$. Again, as $D$ is a dominating set we must have $p_1z \in E(G)$, and $p_kz \in E(G)$. Consider the set $\{x,y, p_1, p_k, w\}$. Since $w\notin V(P)$ neither of $wp_1$ nor $wp_k$ is an edge, similarly, since $x,y \notin V(P)$, $\{x,y,p_1,p_k\}$ induces a $P_2+2P_1$. Moreover, by Observation~\ref{obs:ifp1pkinE}, $p_1p_k \in E(G)$. Hence, $\{x,y, p_1, p_k, w\}$ induces a $P_2 + 3P_1$ in $G$, a contradiction.
The case for $V(P) \cap V(D) = \{w\}$ is symmetrical. 
This finishes the proof of \Cref{clm:intersection-geq-2}.
\end{proof}

\begin{claim}\label{clm:p-intersect-xy-isnotempy}
    $V(P) \cap \{x,y\} \neq \emptyset$.
\end{claim}

\begin{proof}[Proof of \Cref{clm:p-intersect-xy-isnotempy}]\renewcommand{\qedsymbol}{}
Suppose not. Let $y'\in V(G)$ be a vertex in $V(P)$ that is closest to $\{x,y\}$ and let $Q$ be a shortest path between $\{x,y\}$ and $y'$ (as $G$ is connected). We may assume that $y$ and $y'$ are the endpoints of this path. Let $\hat{y}$ be the vertex adjacent to $y'$ in $Q$. Note that it may be that $\hat{y}=y$. Let $c$ and $d$ be the vertices before and after $y'$ in $P$. By \Cref{clm:intersection-geq-2}, we have $z\in V(P)$ and $w \in V(P)$. Let $a$ and $b$ be the neighbors of $z$ in $P$. Observe that neither $p_1x \in E(G)$ nor $p_1y \in E(G)$, as otherwise $P$ could be extended to $x$ or $y$. Analogously, $p_kx \notin E(G)$ and $p_ky \notin E(G)$. Since $D$ is a dominating set, $p_1$ and $p_k$ must be dominated by $\{z,w\}$. Therefore we have two cases.

\begin{figure}[t]
    \centering
    \includegraphics[scale=1.3]{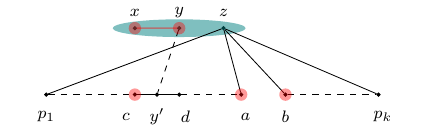}
    \caption{Case~\ref{case:p1pk-dominatedbysame-vrtx}: $p_1,p_k$ dominated by $z$. Dashed lines show paths of arbitrary length.}
   \label{fig:case1}
\end{figure}

\begin{case}\label{case:p1pk-dominatedbysame-vrtx}
    There exists a vertex $q\in\{z,w\}$ such that both $p_1$ and $p_k$ are dominated by $q$.
\end{case}

Suppose, without loss of generality, that $q=z$ and $z$ is between $y'$ and $p_k$ in $P$. Since $p_1$ and $p_k$ dominated by $z$, we may assume that $ab \notin E(G)$, otherwise the path $xy\dots y'\dots p_1zp_k \dots ba \dots d$ is a longer path. We will show that $\{x,y, a, b, c\}$ induces a $P_2 + 3P_1$. First, if $xa\in E(G)$, then $yxa\dots p_1zb\dots p_k$ is a longer path (symmetrically, $bx \notin E(G)$) and if $ya \in E(G)$, then $xya\dots p_1zb \dots p_k$ is a longer path (symmetrically, $by \notin E(G)$). Second, $ac, bc \notin E(G)$, since otherwise $xy\dots y' \dots ac\dots p_1zb \dots p_k$, and $xy\dots y' \dots azp_k \dots bc \dots p_1$ are longer paths.
Third, $cx \notin E(G)$, since otherwise $p_1\dots cxy\dots y'\dots p_k$ is a longer path. We may also assume that $cy \notin E(G)$, as $y\notin V(P)$. Hence, the fact that $ab\notin E(G)$ implies that $\{x,y, a, b, c\}$ form a $P_2 + 3P_1$ (red vertices in \Cref{fig:case1}). This yields a contradiction and shows that \Cref{clm:p-intersect-xy-isnotempy} holds for \Cref{case:p1pk-dominatedbysame-vrtx}.

\begin{case}\label{case:p1dombyz-pkdombyw}
$p_1$ and $p_k$ are dominated by distinct vertices among $\{z,w\}$.  
\end{case}

\noindent Suppose, without loss of generality, that $p_1$ is dominated by $z$ and $p_k$ is dominated by $w$.
We have two possibilities:

\begin{enumerate}[({\rm}a)]
    \item\label{case2a} $y'$ is between $p_1$ and $z$ (symmetrically, $y'$ is between $p_k$ and $w$),    
    \item\label{case2b} $y'$ is between $z$ and $w$.
\end{enumerate}
\bigskip

To prove \Cref{case:p1dombyz-pkdombyw}\ref{case2a} we show that (see \Cref{fig:case2a}):

\clm{$p_1a \in E(G)$.}\label{sta:p_1a}
\begin{proof}[Proof of (\ref{sta:p_1a})]\renewcommand{\qedsymbol}{}
Recall that $a,b\in V(G)$ are the neighbors of $z$ in $P$, and that $\hat{y}$ is the vertex of $Q$ adjacent to $y'$ (and it can be that $\hat{y}=y$). We show that $\{\hat{y},y', p_1, p_k, a\}$ induces a $P_2 + 3P_1$, unless $p_1a \in E(G)$. It is easy to see that, $cd\in E(G)$ by Observation~\ref{obs:bridge1}, and $\hat{y}p_1,\hat{y}p_k \notin E(G)$ as $\hat{y}\notin V(P)$. First, we may assume that $p_1y' \notin E(G)$ otherwise, $xy\dots y'p_1\dots cd \dots p_k$ is a longer path. Similarly, $p_ky' \notin E(G)$ otherwise, $xy\dots y'p_k\dots dc \dots p_1$ is a longer path. Second, neither of $\hat{y}a$ nor $y'a$ is an edge in $E(G)$ since otherwise in the former case, $xy\dots \hat{y}a\dots p_1zb \dots p_k$ is a longer path and in the latter case, $xy\dots y'a\dots dc \dots p_1zb \dots p_k$ is a longer path. From Observation~\ref{obs:ifp1pkinE} we have $p_1p_k \notin E(G)$. Moreover, $p_ka \notin E(G)$ otherwise, $xy\dots y' \dots ap_k \dots bzp_1 \dots c$ is a longer path. This shows that $\{\hat{y},y', p_1, p_k, a\}$ induces a $P_2 + 3P_1$ which is a contradiction, hence $p_1a \in E(G)$. This proves (\ref{sta:p_1a}).
\end{proof}

Now we claim that $\{x,y,z,c,p_k\}$ induces a $P_2 + 3P_1$ in $G$. Observe that $xc \notin E(G)$ since otherwise $p_1\dots cxy\dots y'\dots p_k$ is a longer path and $yc \notin E(G)$ otherwise, we can extend $P$ to $y$. Moreover, since $ap_1 \in E(G)$ by (\ref{sta:p_1a}), neither of $cz$ nor $cp_k$ is an edge in $E(G)$ since otherwise in the former case, $xy\dots y' \dots ap_1 \dots cz\dots p_k$ is a longer path, and in the latter case, $xy\dots y' \dots ap_1 \dots cp_k \dots bz$ is a longer path. Moreover, by our assumption, neither of $xp_k$ nor $yp_k$ is an edge in $E(G)$. If $zp_k \in E(G)$ then $p_1,p_k$ are dominated by $z$ and we are done by \Cref{case:p1pk-dominatedbysame-vrtx}. This implies that $\{x,y, z, c, p_k\}$ induces a $P_2 + 3P_1$ (red vertices in \Cref{fig:case2a}) which is a contradiction. Therefore, Claim~\ref{clm:p-intersect-xy-isnotempy} holds for \Cref{case:p1dombyz-pkdombyw}\ref{case2a}. We leave it to the reader to check that the case $y'$ is between $p_k$ and $w$ can be shown symmetrically.

\medskip

We proceed with \Cref{case:p1dombyz-pkdombyw}\ref{case2b}, that is, when $y'$ is between $z$ and $w$. We remind the reader that by Lemma~\ref{lem:three-cut-vertex}, we assume $G$ has no set of three vertices that are cut vertices and induce a triangle. Our goal is to show that, in this case, such vertices exist, leading to a contradiction. Note that, $w\in V(P)$ and $cd\in E(G)$ by Observation~\ref{obs:bridge1}. Let $e,f \in V(G)$ be the neighbors of $w$ in $P$.
Let $c'$ be the vertex before $c$, and $d'$ be the vertex after $d$ in $P$ (see Figure~\ref{fig:adj}).

\begin{figure}[t]
    \centering
    \includegraphics[scale=1.3]{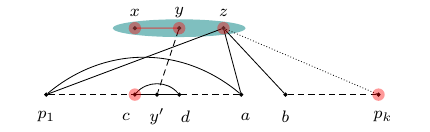}
    \caption{\Cref{case:p1dombyz-pkdombyw}\ref{case2a}: $y'$ is between $p_1$ and $z$. Dashed lines show paths of arbitrary length, dotted lines show non-edges.}
   \label{fig:case2a}
\end{figure}

\begin{figure}[t]
    \centering
    \includegraphics[scale=1.5]{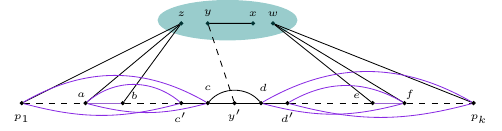}
    \caption{\Cref{case:p1dombyz-pkdombyw}\ref{case2b}: $y'$ is between $z$ and $w$, outcomes of (\ref{adj}). Dashed lines show paths of arbitrary length.}
    \label{fig:adj}
\end{figure}

\clm{The edges $ac, df, ac', p_1c', p_kd, p_kd', p_1c, d'f, ab,$ and $ef$ exist in $G$.}\label{adj}

\begin{proof}[Proof of (\ref{adj})]\renewcommand{\qedsymbol}{}
Assume for the sake of contradiction that it is not the case, then:

    \begin{enumerate}[({\rm}I)]\setlength\itemsep{0.5em}
    
        \item\label{adj_E1} if $ac \notin E(G)$, we will show that $\{x,y,a,f,c\}$ induces a $P_2+3P_1$. First, since $x\notin V(P)$ we have $ax\notin E(G)$ otherwise, $yxa\dots p_1zb \dots p_k$ is a longer path. Similarly, $ay\notin E(G)$ otherwise, $xya\dots p_1zb \dots p_k$ is a longer path. Second, $xf\notin E(G)$ otherwise, $yxf\dots p_kwe\dots p_1$ is a longer path, and similarly $yf\notin E(G)$ otherwise, $xyf\dots p_kwe\dots p_1$ is a longer path. Observe that neither of $xc$ nor $yc$ is an edge otherwise, one can include $x$ and $y$ (in the former) or $y$ (in the latter) to $P$. Third, if $af \in E(G)$ then $xy\dots y'\dots bzp_1\dots af\dots p_kwe\dots d$ is a longer path, and if $cf\in E(G)$ then $xy\dots y'\dots ewp_k\dots fc\dots p_1$ is a longer path. This implies that $\{x,y, a,f,c\}$ induces a $P_2+3P_1$ which is a contradiction.
        
        \item\label{adj_E2} If $df \notin E(G)$, take the set $\{x,y, a,f,d\}$. A symmetric argument to \ref{adj_E1} shows that $\{x,y, a,f,d\}$ induces a $P_2+3P_1$ which is a contradiction.
    
        \item\label{adj_E3} If $ac' \notin E(G)$, take the set $\{x,y, a,f,c'\}$. By \ref{adj_E2} it is enough to show that the edges $xc'$, $yc'$, $fc'$ do not exist. Since from \ref{adj_E1}, $ac \in E(G)$, if $xc' \in E(G)$ then $yxc'\dots bzp_1\dots ac\dots p_k$ is a longer path. Similarly, $yc' \notin E(G)$. Moreover $fc' \notin E(G)$ otherwise, $xyy'cd\dots ewp_k\dots fc'\dots p_1$ is a longer path. This shows that $\{x,y, a,f,c'\}$ induces a $P_2+3P_1$ which is a contradiction.

        \item\label{adj_E4} If $d'f \notin E(G)$ then consider the set $\{x,y,a,f,d'\}$. A symmetric argument to \ref{adj_E3} shows that $\{x,y,a,f,d'\}$ induces a $P_2+3P_1$, a contradiction.
        
        \item\label{adj_E5} Let $p_1c' \notin E(G)$ and take the set $\{y,y', c',d',p_1\}$. As we argued, we may assume that $\hat{y}=y$ and so $yy'\in E(G)$. First, as we have shown in \ref{adj_E3} $yc' \notin E(G)$, and $yp_1 \notin E(G)$ as $y\notin V(P)$. Second, if $yd'\in E(G)$ then since $df \in E(G)$ by \ref{adj_E2} we have $xyd'\dots ewp_k \dots fd\dots p_1$ is a longer path. Third, if $y'c' \in E(G)$ then since $ac\in E(G)$ by \ref{adj_E1} we have $xyy'c'\dots bzp_1\dots acd \dots p_1$ is a longer path. If $y'd' \in E(G)$ then $xyy'd' \dots ewp_k \dots fdc \dots p_1$ is a longer path. And if $y'p_1 \in E(G)$ then $xyy'p_1\dots cd\dots p_k$ is a longer path. Finally, neither of $c'd'$ nor $d'p_1$ is an edge otherwise in the former case, as $ca$ is an edge by \ref{adj_E1}, we have $xyy'dca\dots p_1zb \dots c'd'\dots p_k$ is a longer path and in the latter case, $xyy'dc\dots p_1d'\dots p_k$ is a longer path. This shows that $\{y,y', c',d',p_1\}$ induces a $P_2+3P_1$ which is a contradiction.

        \item\label{adj_E6} If $p_kd' \notin E(G)$, take the set $\{p_1,c', x,d',p_k\}$. By a symmetric argument to \ref{adj_E5} one can observe that $\{p_1,c', x,d',p_k\}$ induces a $P_2+3P_1$, a contradiction.

        \item\label{adj_E7} Let $p_kd \notin E(G)$ and consider the set $\{p_1,c', x,d,p_k\}$. Observe that by \ref{adj_E5} we have $p_1c' \in E(G)$. First, $p_1x \notin E(G)$ as $x \notin V(P)$ and $p_1p_k \notin E(G)$ by Observation~\ref{obs:ifp1pkinE}. If $p_1d\in E(G)$, then $xyy'\dots p_1d \dots p_k$ is a longer path. Second, as we have shown in \ref{adj_E3}, $c'x \notin E(G)$, and $c'd \notin E(G)$ otherwise, as $ca$ is an edge by \ref{adj_E1}, we have $xyy'ca\dots p_1zb\dots c'd \dots p_k$ is a longer path. Similarly, if $c'p_k \in E(G)$, then $xyy'ca\dots p_1zb\dots c'p_k \dots d$ is a longer path. Finally, observe that neither of $xd$ nor $xp_k$ is an edge in $G$ as otherwise, in the former case, $p_k\dots dxyy' \dots p_1$ is a longer path, the latter case also follows by the fact that $x \notin V(P)$. This along with the fact that $p_1c' \in E(G)$ proved in \ref{adj_E5}, implies that $\{p_1,c', x,d,p_k\}$ induces a $P_2+3P_1$ which is a contradiction.

      \item\label{adj_E8} Let $p_1c \notin E(G)$ and consider the set $\{p_k,d', x,c,p_1\}$. By a symmetric argument to \ref{adj_E7} one can observe that $\{p_k,d', x,c,p_1\}$ induces a $P_2+3P_1$, a contradiction.
    
        \item\label{adj_E9} Suppose $ab \notin E(G)$ and take the set $\{p_k,d', x, a, b\}$. Note that by \ref{adj_E6} we have $p_kd' \in E(G)$. First, since $x\notin V(P)$ we have $p_kx\notin E(G)$. Moreover, $p_ka \notin E(G)$, otherwise $xyy'\dots bzp_1 \dots ap_k \dots d$ is a longer path. Also $p_kb \notin E(G)$, otherwise since $ca\in E(G)$ by \ref{adj_E1}, $xyy'\dots p_kb \dots ca\dots p_1z$ is a longer path. Second, if $d'x \in E(G)$ then as $p_kd$ is an edge by \ref{adj_E7} we have $yxd' \dots p_kd \dots p_1$ is a longer path which is impossible, hence $d'x \notin E(G)$. If $ad' \in E(G)$ then $xyy' \dots bzp_1 \dots ad' \dots p_kd$ is a longer path, hence $d'a \notin E(G)$. If $d'b \in E(G)$ then it follows by \ref{adj_E7} and \ref{adj_E8} that $dp_k,cp_1 \in E(G)$, and $xyy'dp_k \dots d'b \dots cp_1 \dots az$ is a longer path. Third, as we have shown in \ref{adj_E1}, $xa \notin E(G)$. Moreover, since $cp_1$ is an edge, if $xb \in E(G)$ then $p_k\dots y'yxb\dots cp_1 \dots az$ is a longer path which is impossible. Hence, $\{p_k, d', x, a, b\}$ induces a $P_2+3P_1$ which is a contradiction. This implies that $ab \in E(G)$.

        \item\label{adj_E10} Finally, suppose $ef \notin E(G)$ and take the set $\{p_1,c', x, e, f\}$. By a symmetric argument to \ref{adj_E9} one can show that $\{p_1,c', x, e, f\}$ induces a $P_2+3P_1$ which is a contradiction.
\end{enumerate}

    This, completes the proof of (\ref{adj}).
\end{proof}

For the rest of the proof of \Cref{case:p1dombyz-pkdombyw}\ref{case2b}, we wish to show that $\{c,y',d\}$ are three cut vertices that induce a triangle in $G$ (see \Cref{fig:adj}). Let $w' \in V(G)$ be any vertex from $p_1$ to $y'$, and $z' \in V(G)$ be any vertex from $y'$ to $p_k$. We say $\ell$ is a \textit{crossing edge} in $G$ if either $\ell=zz'$ or $\ell=ww'$.

\clm{$G$ has no crossing edges.}\label{sta:crossing-edges}

\begin{proof}[Proof of (\ref{sta:crossing-edges})]\renewcommand{\qedsymbol}{}
Assume for the sake of contradiction that $z$ has a neighbor $z'$ from $y'$ to $p_k$. That is, $\ell = \{z,z'\}$ be a crossing edge in $G$. Let $z^{+}$ be a vertex after $z'$ in $P$. First, note that $z' \neq p_k$ as we are in \Cref{case:p1dombyz-pkdombyw}, and $z' \neq y'$ as otherwise $xyy'zp_1 \dots cd \dots p_k$ is a longer path. Note that $ab \in E(G)$ by \ref{adj_E9}, and $dp_k \in E(G)$ by \ref{adj_E5}. Then it follows that the path $xyy' \dots ba \dots p_1zz' \dots dp_k \dots z^{+}$ is longer than $P$. Analogously, let $w$ has a neighbor $w'$ from $p_1$ to $y'$ and $w^{-}$ be a vertex before $w$ in $P$. Again, as we are in \Cref{case:p1dombyz-pkdombyw}, $w' \neq p_1$, and $w' \neq y'$ as otherwise, $xyy'wp_k\dots dc \dots p_1$ is a longer path which is impossible. Note that $ef \in E(G)$ by \ref{adj_E10}, and $p_1c \in E(G)$ by \ref{adj_E7}.

Then it follows that the path $xyy' \dots ef \dots p_kww'\dots cp_1\dots w^{-}$ is longer than $P$. This is a contradiction, hence $G$ has no crossing edges.
\end{proof}

\clm{$x$ and $y$ have no neighbor in $P$ other than $y'$.}\label{xhasnoneihbor}

\begin{proof}[Proof of (\ref{xhasnoneihbor})]\renewcommand{\qedsymbol}{}
Suppose not and let $x'$ be a neighbor of $x$ in $P$. Recall that $x'\neq p_1$ and $x' \neq c$, otherwise we can extend $P$ to $x$. If $x'$ is between $p_1$ and $c$, since $cp_1$ is an edge by \ref{adj_E7}, $p_k \dots y'yxx' \dots cp_1 \dots x'$ is a longer path. Analogously, let $x'$ be between $d$ and $p_k$ and see that, $x'\neq d$ and $x'\neq p_k$, then as $dp_k$ is an edge by \ref{adj_E5}, $p_1\dots y'yxx' \dots dp_k \dots x'$ is a longer path which is impossible. Now let $y''\neq y'$ be a neighbor of $y$ in $P$. Since $y\notin V(P)$ we have $y''\neq p_1, p_k$, and $y''\neq c,d$, otherwise we can extend $P$ to $y$. Hence, we may assume that $y''$ is a vertex either between $p_1$ and $c$ or between $d$ and $p_k$. In the former case, $p_k \dots y'yy''\dots cp_1\dots y''$, and in the latter case, $p_1\dots y'yy''\dots dp_k\dots y''$ are longer paths. This is a contradiction, hence the claim of (\ref{xhasnoneihbor}) holds.
\end{proof}

\clm{$y'$ has no neighbor from $p_1$ to $c'$, and from $p_k$ to $d'$.}\label{y'hasnoneighbor}

\begin{proof}[Proof of (\ref{y'hasnoneighbor})]\renewcommand{\qedsymbol}{}
Assume for the sake of contradiction that $y'$ has a neighbor $y''$ in $P$. First, observe that neither $y'' = c'$ nor $y'' = p_1$ since otherwise as $cp_1$ is an edge by \ref{adj_E7}; in the former case, $p_k \dots dcp_1 \dots c'y'yx$ and in the latter case, $p_k \dots dc\dots p_1y'yx$ is a longer path. Similarly, $y''\neq d'$ and $y''\neq p_k$.
Thus we may assume that $y''$ is a vertex between $p_1$ and $c'$, or between $d'$ and $p_k$. We suppose that $y''$ is a vertex between $p_1$ and $c'$ and we leave it to the reader to check that the case $y''$ is between $d'$ and $p_k$ is symmetric.
Let $y^{+}, y^{-}$ be the vertices after and before $y''$ in $P$. We show that:

\begin{enumerate}[({\rm}A)]\setlength\itemsep{0.6em}
    \item\label{y'-hasno-neighbor-first} $y''c \in E(G)$: suppose not, then from (\ref{xhasnoneihbor}) we have $yc \notin E(G)$, and $yy'' \notin E(G)$, hence $\{y,y',y'',c \}$ induces a claw which is a contradiction.
    \item\label{y'-hasno-neighbor-second} $y''d \in E(G)$: suppose not. Again from (\ref{xhasnoneihbor}) we have $yd \notin E(G)$, $yy'' \notin E(G)$, hence $\{y,y',y'',d \}$ induces a claw which is a contradiction.
    \item\label{y'-hasno-neighbor-third} $y^{+}y^{-}\in E(G)$: suppose not, we show that $\{y^{-},y'',y^{+},y'\}$ induces a claw. If $y^{-}y' \in E(G)$ then from \ref{y'-hasno-neighbor-first} we have, $xyy'y^{-}\dots p_1c' \dots y''cd \dots p_k$ is a longer path. Analogously, if $y^{+}y'\in E(G)$ then, $xyy'y^{+}\dots c'p_1 \dots y''cd \dots p_k$ is a longer path, hence $\{y^{-},y'',y^{+},y'\}$ induces a claw which is impossible. This implies that $y^{+}y^{-}\in E(G)$.
\end{enumerate}

\medskip

Since $y^{-}y^{+}\in E(G)$ by \ref{y'-hasno-neighbor-third}, we will show that $\{y^{-},y^{+},c,x,p_k\}$ induces $P_2+3P_1$ in $G$. We may assume that, by (\ref{xhasnoneihbor}), that the edges $y^{-}x, y^{+}x, cx, xp_k$ do not exist. If $cp_k\in E(G)$ then $xyy' \dots p_kc \dots p_1$ is a longer path. Also $\{y^{+},y^{-},p_k\}$ induces a $P_2 + P_1$ in $G$, otherwise if $y^{+}p_k \in E(G)$ then, by \ref{y'-hasno-neighbor-first}, $xyy' \dots p_ky^{+} \dots cy'' \dots p_1$ is a longer path and if $y^{-}p_k \in E(G)$ then $xyy' \dots p_ky^{-} \dots cp_1 \dots y^{-}$ is a longer path (see that $cp_1\in E(G)$ by \ref{adj_E7}). It remains to discuss that neither $y^{-}c$ nor $y^{+}c$ is an edge, otherwise in the former case, $p_k\dots dcy^{-}\dots p_1c' \dots y''y'yx$ is a longer path (see that $p_1c'\in E(G)$ by \ref{adj_E4}) and in the latter case, $p_k \dots dcy^{+}\dots c'p_1 \dots y''y'yx$ is a longer path which is impossible. This shows that $\{y^{-},y^{+},c,x,p_k\}$ induces a $P_2+3P_1$, a contradiction. Hence, the claim of (\ref{y'hasnoneighbor}) holds.
\end{proof}

\clm{There is no edge $e=\alpha\beta$ with $\alpha$ from $p_1$ to $c'$, and $\beta$ from $d'$ to $p_k$. }\label{noedgebetween-left-right}

\begin{proof}[Proof of (\ref{noedgebetween-left-right})]\renewcommand{\qedsymbol}{}
Assume for the sake of contradiction that $\alpha\beta \in E(G)$. As before, we denote by $\alpha^{-}, \alpha^{+}$ the vertices before and after $\alpha$, and by $\beta^{-}, \beta^{+}$ the vertices before and after $\beta$ in $P$. We may assume that $\alpha \neq p_1, c'$, otherwise since $dp_k\in E(G)$ by \ref{adj_E5} and $cp_1\in E(G)$ by \ref{adj_E7} we have, $xyy' \dots p_1\beta \dots dp_k \dots \beta^{+}$ and $xyy'cp_1 \dots c'\beta \dots dp_k\beta^{+}$ are longer paths (symmetrically, $\beta \neq p_k,d'$). We wish to show that $\{\alpha^{+},p_1,p_k,x,y\}$ form a $P_2+3P_1$ in $G$. First, observe that $p_1p_k \notin E(G)$ by Observation~\ref{obs:ifp1pkinE}, and by (\ref{xhasnoneihbor}) we have $\{x,y,p_1\}$ and $\{x,y,\alpha^{+}\}$ induces a $P_2+P_1$. Second, neither of $\alpha^{+}p_k$ nor $\alpha^{+}p_1$ is an edge, otherwise in the former case, $xyy'\dots p_k\alpha^{+}\dots cp_1 \dots \alpha$ is a longer path and in the latter case, $xyy' \dots \alpha^{+}p_1\dots \alpha\beta \dots dp_k \dots \beta^{+}$ is a longer path. This shows that $\{\alpha^{+},p_1,p_k,x,y\}$ induces a $P_2+3P_1$ which is a contradiction. In a symmetrical way one could show that $\{\beta^{-},p_1,p_k,x,y\}$ induces a $P_2+3P_1$. This completes the proof of (\ref{noedgebetween-left-right}).
\end{proof}

It is now easy to conclude from (\ref{sta:crossing-edges})-(\ref{noedgebetween-left-right}) that the vertices in $K_3 = cy'd$ are cut vertices in $G$. This violates the assumption that $G$ has no three cut vertices inducing a triangle, hence completing the proof of Claim~\ref{clm:p-intersect-xy-isnotempy}.
\end{proof}

We are ready to finish the proof of Theorem~\ref{thm:claw,p2+3p1}. Recall that our desire is to show that $x$ is a Gallai vertex. Since $x \notin V(P)$ (by the assumption), Claim~\ref{clm:p-intersect-xy-isnotempy} implies that $y\in V(P)$, say $y=p_i$ for $1<i<k$. As $d_{G}(x)\geq d_{G}(y)$, $x$ must have a private neighbor in $G$, say $x'$. Recall that, by Observation~\ref{obs:bridge1} we have $p_{i-1}p_{i+1}\in E(G)$ as $x\notin V(P)$. We have two possibilities:

\begin{itemize}\setlength\itemsep{0.5em}
    \item If $x' \notin V(P)$ then we show that $\{p_1,p_2,x',p_k,y\}$ induces a $P_2+3P_1$. First, neither of $p_1x'$ nor $p_2x'$ is an edge; otherwise, in the former case $yxx'p_1\dots p_{i-1}p_{i+1}\dots p_k$ is a longer path, and in the latter case $yxx'p_2 \dots p_{i-1}p_{i+1} \dots p_k$ is a longer path, symmetrically, $p_kx'\notin E(G)$. Observe that since $x'$ is a private neighbor of $x$, $yx' \notin E(G)$. Second, if $p_2p_k \in E(G)$ then the path $x'xyp_{i-1}\dots p_2p_k\dots p_{i+1}$ is a longer path. Third, neither of $p_1y$ nor $p_2y$ is an edge since otherwise, in the former case $x'xyp_1\dots p_{i-1}p_{i+1} \dots p_k$ is a longer path, and in the latter case the path $x'xyp_2\dots p_{i-1}p_{i+1} \dots p_k$ is a longer path, symmetrically, $p_ky\notin E(G)$. This shows that $\{p_1,p_2,x',p_k,y\}$ induces a $P_2+3P_1$ in $G$ which is impossible.

\item By the first bullet point we may suppose that $x' \in V(P)$. Let $x'$ be a vertex between $p_{i+1}$ and $p_k$, and observe that the case for when $x'$ is between $p_{i}$ and $p_1$ is symmetric. We denote by $x^{-},x^{+}$ the vertices before and after $x'$ in $P$. We show that $\{x,y,p_1,p_k,x^{-}\}$ induces a $P_{2}+3P_{1}$ in $G$. Since $G$ is claw-free we have $x^{-}x^{+} \in E(G)$ otherwise, $\{x^{-},x',x^{+},x\}$ induces a claw. Moreover, it follows from Observation~\ref{obs:bridge1} that $\{x,y, x^{-},x^{+}\}$ induces a $2P_2$ in $G$. Note that, $\{x,y,p_1\}$ and $\{x,y,p_k\}$ induces a $P_2+P_1$ in $G$ since otherwise we could extend $P$ to $x$ ($p_{i-1}p_{i+1}\in E(G)$ by Observation~\ref{obs:bridge1}). If $p_1x^{-}\in E(G)$ then $p_{i+1}\dots x^{-}p_1 \dots p_{i-1}yxx' \dots p_k$ is a longer path, and if $p_kx^{-} \in E(G)$ then $p_{i+2}\dots x^{-}p_k \dots x'xyp_{i+1}p_{i-1} \dots p_{1}$ is a longer path. This shows that $\{x,y,p_1,p_k,x^{-}\}$ induces a $P_{2}+3P_{1}$ which is a contradiction.
\end{itemize}

The above two bullet points conclude that $x$ must be a common vertex to all longest paths in $G$. This finishes the proof of Theorem~\ref{thm:claw,p2+3p1}.
\end{proof} 

\begin{observation}\label{alg:claw,P2+3P1}
Let $G$ be a connected $(\claw, P_2+3P_1)$-free graph. There exists a polynomial-time algorithm which finds a Gallai vertex in $G$.
\end{observation}

The proof of Observation~\ref{alg:claw,P2+3P1} is quite similar to the proof of Observation~\ref{alg:claw,P3+2P1}. We remind the reader that in a connected $(P_2 + 2P_1)$-free graph, every vertex of maximum degree is a Gallai vertex (see Proposition 7 in \cite{long2023non}).

\medskip

Let us now put everything together. Theorems~\ref{thm:claw,p3+2p1}, \ref{thm:claw,K3+2P1}, \ref{thm:claw,2P_2+P_1}, along with Theorem~\ref{thm:claw,p2+3p1} yields Theorem~\ref{thm:main-claw,Hinintgro}.

\begin{figure}[t]
    \centering
    \includegraphics[scale=1.1]{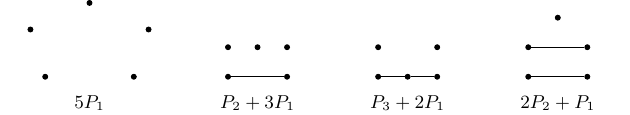}
    
           \bigskip
           
    \includegraphics[scale=1]{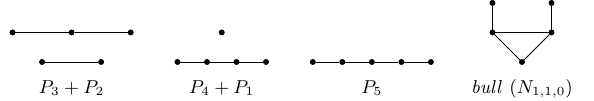}
    
         \bigskip
         
    \includegraphics[scale=1]{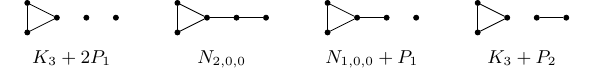}

    \caption{All five-vertex induced subgraphs of $\mathcal{B}^{+}$.}
    \label{fig:5-vertexgraphs}
\end{figure}

Now Theorem~\ref{thm:main-claw,Hinintgro}, combined with Theorem~\ref{thm:dich-for-claw-H} and Theorem~\ref{2p2-free} implies a complete classification, for all graphs $H$ of size at most five, of $(\claw, H)$-free graphs that form a Gallai family. More precisely:

\begin{corollary}\label{corr:H,claw}
Let $H$ be a graph on at most five vertices. The class of $(claw,H)$-free graphs forms a Gallai family if and only if one of the following holds.

    \begin{itemize}
        \item $H$ is a linear forest;
        \item $H$ is one of the graphs $bull, K_3+2P_1, N_{2,0,0}, N_{1,0,0}+P_1$ or $K_3+P_2$.
    \end{itemize}
\end{corollary}

\begin{proof}
For the first bullet point, necessity follows from Theorem~\ref{necessityforlinearforests}. For the sufficiency part, we wish to show that if $H$ is a linear forest of size at most five then $(\claw, H)$-free graphs form a Gallai family. We have the following cases: $5P_1$, $P_5$, $P_2+3P_1$, $P_3+2P_1$, $2P_2+P_1$, $P_3+P_2$, and $P_4+P_1$ (see \Cref{fig:5-vertexgraphs}). Observe that the case $H=5P_1$ has been settled in~\cite{long2023non}. Moreover, when $H=P_5$, $H=P_3+P_2$ and $H=P_4+P_1$, then $(claw,H)$-free graphs are a subclass of $(claw,P_6)$-free graphs which were shown by Gao and Shan~\cite{gao2021nonempty} to form a Gallai family (see Theorem~\ref{thm:dich-for-claw-H}). The remaining cases, $P_2+3P_1$, $P_3+2P_1$, and $2P_2+P_1$ conclude from Theorem~\ref{thm:main-claw,Hinintgro}. 

For the second bullet point, to see the necessity part let $|V(H)|=5$ and suppose that $H$ is not a linear forest.
Since $\mathcal{B}^{+}$ is claw-free we may assume that $H$ is an induced subgraph of $\mathcal{B}^{+}$; otherwise, it violates the assumption that the class of $(claw, H)$-free graphs form a Gallai family. It is easy to observe that (\Cref{fig:5-vertexgraphs}) $bull, K_3+2P_1, N_{2,0,0}, N_{1,0,0}+P_1$ and $K_3+P_2$ are the only five-vertex (except linear forests) induced subgraphs of $\mathcal{B}^{+}$. For the sufficiency part, observe that the case $H= bull$ and $H=N_{2,0,0}$ has been settled in Theorem~\ref{thm:dich-for-claw-H}. Note that, when $H=K_3+P_2$ and $H=N_{1,0,0}+P_1$ then $(claw,H)$-free graphs are a subclass of $(claw,N_{3,0,0})$-free graphs. The case $H=K_3+2P_1$ follows directly from Theorem~\ref{thm:main-claw,Hinintgro}. This completes the proof of Corollary~\ref{corr:H,claw}.
\end{proof}

\section{Gallai family in $(P_5,H)$-free graphs}\label{sec:p5,H}
We begin the section with proving Theorem~\ref{thm:main} for $H=K_3$. The reader may observe that the class of $(P_5, \diam)$-free graphs contains the class of $(P_5, \trig)$-free graphs. However, we include a proof for $(P_5, \trig)$-free case as it is
a vital piece in the proof of $(P_5, \paw)$-free case.

\begin{lemma}[\cite{bacso1990dominating}]\label{lem:p5,c5}
A graph $G$ is $(P_5,C_5)$-free if and only if every induced subgraph of $G$ contains a dominating clique.
\end{lemma}

A 5-$ring$ is a graph whose vertices can be partitioned into five non-empty stable sets $S_1, \dots, S_5$ such that for each $i$ (in modulo 5) $S_i$ is complete to $S_{i-1} \cup S_{i+1}$ and anticomplete to $S_{i-2} \cup S_{i+2}$. Sumner~\cite{sumner1981subtrees} characterized the structure of $(P_5, \text{triangle})$-free graphs. For completeness, we include a short proof below. 

\begin{lemma}\label{lem:c5-structure}
Let $G$ be a $(P_5,\text{triangle})$-free graph. Then each component of $G$ is either a 5-ring or bipartite.
\end{lemma}

\begin{proof}[Proof of \Cref{lem:c5-structure}]
Suppose not, that is, let $H$ be a component of $G$ that contains an induced $C_5$. So we may assume that $H$ contains a 5-ring, say $C= S_1, \dots, S_5$. Suppose that $C$ has the maximum size in $H$. Let $H \neq C$, and pick a vertex $s_i$ in each $S_i$, $1\leq i\leq 5$. Let $w$ be a vertex in $H\setminus C$ with a neighbor $s_1 \in S_1$. Since $G$ is triangle-free neither of $ws_2$ nor $ws_5$ is an edge in $G$. So $w$ has no neighbor in $S_2 \cup S_5$. Also if $ws_4 \notin E(G)$ (symmetrically $ws_3 \notin E(G)$) then $\{w,s_1,s_2,s_3,s_4\}$ induces a $P_5$, therefore we may assume that $w$ is complete to $S_1 \cup S_4$ (symmetrically to $S_1 \cup S_3$). But then $S_1,\dots,S_5\cup \{w \}$ (symmetrically $S_1, S_2\cup \{w\}, \dots, S_5$) is a 5-ring which is a contradiction. Thus, $H = C$ and $H$ is a 5-ring. This proves \Cref{lem:c5-structure}.
\end{proof}

\medskip

Note that a 5-ring graph is $2P_2$-free. We are now ready to prove the following:

\begin{theorem}\label{thm:p5,tring}
The class of $(P_5, \trig)$-free graphs forms a Gallai family.
\end{theorem}

\begin{proof}
Let $G$ be a connected $(P_5, \trig)$-free graph. If $G$ contains an induced $C_5$, then it follows from Lemma~\ref{lem:c5-structure} that $G$ is a 5-ring, hence $2P_2$-free, and Theorem~\ref{2p2-free} implies that $G$ has a Gallai vertex. Therefore, we may assume that $G$ is a $(P_5, C_5, \trig)$-free graph. Then by Lemma~\ref{lem:p5,c5}, $G$ contains a dominating clique, 
which has size at most two, since $G$ is a triangle-free. If $G$ has a dominating vertex, by Observation~\ref{obs:domupperbound}, $G$ has a Gallai vertex. Let $D = \{x,y\}$ with $xy\in E(G)$ be a dominating edge in $G$. We may suppose, without loss of generality, that $|N_{G}(y)| \leq |N_{G}(x)|$.
We wish to show that $x$ is a Gallai vertex in $G$. Let $P := p_1p_2 \dots p_k$ be a longest path, and assume for the sake of contradiction that $x$ does not belong to $P$. Since $D$ is a dominating set, $p_1$ must be adjacent to at least one vertex among $\{x,y\}$. Since $P$ is a longest path and $x\notin V(P)$, $p_1x \notin E(G)$ and so $p_1y \in E(G)$. Similarly $p_ky \in E(G)$. We have $y\in V(P)$, otherwise we could extend $P$ to $y$. Let $y=p_i$ for some $1<i<k$. We claim that $|N_{P}(y)| \geq |N_{P}(x)|$. To see this, observe that for each $1<j<k$, either of $p_jx$ or $p_jy$ is an edge. Since $G$ is triangle-free if $p_jx \in E(G)$ then neither $p_{j-1}x\in E(G)$ nor $p_{j+1}x\in E(G)$. This combined with the fact that $y\in V(P)$, and $p_1y, p_ky\in E(G)$ proves our claim.
Now consider the set $\{p_{i-2},p_{i-1},p_{i},p_{i+1},p_{i+2}\}$. Observe that $p_{i-1}p_{i+1}\notin E(G)$ otherwise, $xyp_1\dots p_{i-1}p_{i+1}\dots p_{k}$ is a longer path. Since $G$ is triangle-free, neither of $p_{i-2}p_i$ nor $p_{i+2}p_i$ is an edge. Suppose $p_{i-1}p_{i+2}\in E(G)$. Since $|N_{P}(y)| \geq |N_{P}(x)|$ and $|N_{G}(y)| \leq |N_{G}(x)|$, there is a vertex $u\in N_{G}(x)$ such that $u \notin V(P)$. But then $uxyp_1 \dots p_{i-1}p_{i+2}\dots p_{k}$ is a longer path, as it excludes $p_{i+1}$ and includes $u,x$ in $P$. Analogously, if $p_{i+1}p_{i-2}\in E(G)$ then $uxyp_1 \dots p_{i-2}p_{i+1}\dots p_{k}$ is a longer path, as it excludes $p_{i-1}$ and includes $u,x$ in $P$. Therefore, neither of $p_{i-1}p_{i+2}$ nor $p_{i+1}p_{i-2}$ is an edge. Moreover, observe that since $G$ is $C_5$-free we have $p_{i-2}p_{i+2}\notin E(G)$. This implies that $\{p_{i-2},p_{i-1},p_{i},p_{i+1},p_{i+2}\}$ induces a $P_5$ in $G$ which is a contradiction, hence $x$ is common to all longest paths in $G$, as desired.
 \end{proof}

\begin{observation}\label{alg:P5,K3}
Let $G$ be a connected $(P_5, \trig)$-free graph. There exists a polynomial-time algorithm which finds a Gallai vertex in $G$.
\end{observation}

\begin{proof}[Proof of Observation~\ref{alg:P5,K3}]{}
 As we have shown above, there are two possibilities: either $G$ contains an induced $C_5$ and so $2P_2$-free and then by the main result of \cite{golanshan2k2} every vertex of maximum degree is a Gallai vertex in $G$, or $G$ is a $(P_5, \trig, C_5)$-free graph and so has a dominating clique of size at most two which contains a Gallai vertex. In the former case, it is clear that there is a polynomial-time algorithm to find a Gallai vertex. The latter case follows by a $\mathcal{O}(|V|^3)$-time algorithm to find a dominating clique of size at most two.
\end{proof}

\medskip

Next, we give a short argument for the case $H=\paw$ of Theorem~\ref{thm:main}. Recall that a paw is a $N_{1,0,0}$, the following result describes the structure of paw-free graphs. 

\begin{lemma}[\cite{olariu1988paw}]\label{lem:pawfree}
    A graph $G$ is paw-free if and only if each component of $G$ is triangle-free or complete multipartite.
\end{lemma}

If $G$ is a $(P_5,\paw)$-free graph, then by Lemma~\ref{lem:pawfree}, $G$ is either $(P_5, \trig)$-free or a complete multipartite. In the former case, as we have shown in Theorem~\ref{thm:p5,tring}, $G$ has a Gallai vertex. In the latter case, the fact that multipartite graphs are $2P_2$-free, combined with Theorem~\ref{2p2-free} implies that $G$ has a Gallai vertex. The following is immediate:

\begin{corollary}\label{sta:p5,paw}
The class of $(P_5,\text{paw})$-free graphs forms a Gallai family.
\end{corollary}

\medskip

We now plunge into the $(P_5, \diam)$-free case. We begin by describing three types of graphs that we will need later. Let $\mathcal{G}_{1}, \mathcal{G}_{2}, \mathcal{G}_{3}$ be the families of graphs. For a graph $G$, let

\begin{enumerate}\setlength\itemsep{0.5em}
\item $G \in \mathcal{G}_{1}$ if $G$ contains a dominating clique $K$ such that:
\begin{itemize}
   \item every component of $G\setminus K$ is complete,
    \item for each component $D$ in $G\setminus K$, there is a unique vertex $x\in K$ such that $G[V(D) \cup \{x\}]$ is complete,
    \item the above be the only edges in $G$.
\end{itemize}

\item $G \in \mathcal{G}_{2}$ if $V(G)$ can be partitioned into two sets $X$ and $Y$ with the following specifications:
\begin{itemize}
   \item $G[X]$, $G[Y]$ are complete,
    \item $G[X \cup Y]$ is a non-empty set of independent edges.
\end{itemize}

\item $G \in \mathcal{G}_{3}$ if $G$ contains a $2P_2$-free bipartite graph $H$ with an edge $e=xy$ dominates $G$ such that $G[(V(H)\setminus e), V(G)\setminus V(H)] = \emptyset$, and for all components $D_i$ of $G\setminus V(H)$ the following holds:
\begin{itemize}
    \item for each $i$, $D_i$ is complete,
    \item for each $i$, either $G[V(D_i) \cup \{x\}]$ or $G[V(D_i) \cup \{y\}]$ is complete,
    \item for at most one $i$, $G[V(D_i) \cup \{x\} \cup \{y\}]$ is complete,
    \item the above be the only edges in $G$.
\end{itemize}
\end{enumerate}

\begin{figure}[t]
    \centering
    \includegraphics[scale=1.1]{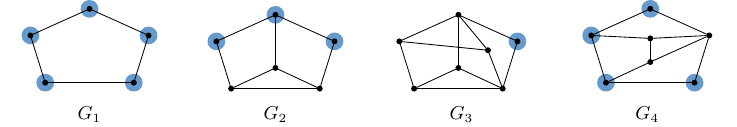}
    \includegraphics[scale=1.1]{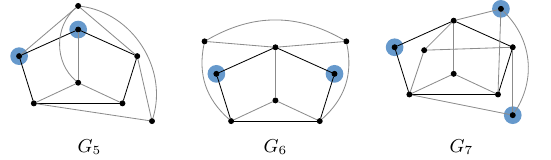}
    \includegraphics[scale=1.1]{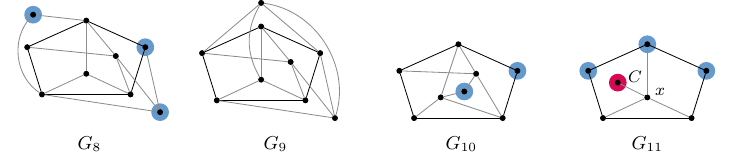}
    \caption{The graphs used in Lemma~\ref{lem:basicgraphs}.}
    \label{fig:BS}
\end{figure}

We need the following decomposition for $(P_5, C_5, \diam)$-free graphs:

\begin{lemma}[\cite{choudum2010first}]\label{lem:p5.c5.diamond}
A graph $G$ is $(P_5, C_5, \diam)$-free if and only if $G \in \mathcal{G}_1 \cup \mathcal{G}_2 \cup \mathcal{G}_3$.
\end{lemma}

We mention two remarks on the families $\mathcal{G}_1,\mathcal{G}_2$. A graph is \textit{star-like} if it is an intersection graph of substars of a star \cite{cerioli2006characterizing}. The following result is due to Cerioli and Lima \cite{cerioli2020intersection}:

\begin{lemma}[\cite{cerioli2020intersection}]
Let $\mathcal{H}$ be a connected graph in which $V(\mathcal{H})$ can be partitioned into $k+1$ sets $K, V_1, \dots, V_k$, for some $k \in \mathbb{N}$, such that the following holds
\begin{itemize}
    \item $K$ is a clique;
    \item For all $x \in V_i$ and $y \in V_j$, $i \neq j$, it holds that $xy \notin E(\mathcal{H})$;
    \item The vertices of $V_i$ can be ordered $v_{i1}, v_{i2}, \dots, v_{i|Vi|}$ in such a way that for all $x \in K$, if $xv_{ij} \in E(\mathcal{H})$, then $xv_{ik} \in E(\mathcal{H})$ for
all $k < j$;
\item For all $x \notin K$, there exists $y \in K$ such that $xy \in E(\mathcal{H})$.
\end{itemize}
Then $\mathcal{H}$ has a Gallai vertex.
\end{lemma}

The reader may observe that the family $\mathcal{G}_{1}$ is contained in the class of star-like graphs. Also it is not hard to see that the family $\mathcal{G}_2$ is Hamiltonian.

We need one more lemma to define our next result precisely.
An \textit{expansion} of a graph $H$ is any graph $G$ such that the vertices of $G$ can be partitioned into $|V(H)|$ nonempty sets $S_{x}$, for $x \in V(H)$, such that $G[S_{x} \cup S_{y}]$ is complete if $xy \in E(H)$, and $G[S_{x} \cup S_{y}] = \emptyset$ if $xy \notin E(H)$. An expansion of a graph is a stable set expansion if each $S_x$ is a stable set, and is a $P_{3}$-free expansion if each $S_x$ induces a $P_3$-free graph. The following characterization for $(P_5, \diam)$-free graphs is due to Arbib and Mosca \cite{arbib2002p5} (see also Theorem 2 in \cite{choudum2010first}):
\begin{lemma}[\cite{arbib2002p5}]\label{lem:basicgraphs}
$G$ is a $(P_5, \diam)$-free graph that contains an induced $C_5$ if and only if $G$ is obtained from one of the graph among $G_1, \dots, G_{11}$ by stable set expansion of each blue vertex, and $P_3$-free expansion of the red vertex.
\end{lemma}

We are now ready to prove the following result:

\begin{theorem}\label{thm:p5,diam}
    The class of $(P_5, \diam)$-free graphs forms a Gallai family.
\end{theorem}

\begin{proof}
Let $G$ be a $(P_5, \diam)$-free graph. We first show that the claim of Theorem~\ref{thm:p5,diam} holds for when $G$ has no induced $C_5$. Let $G$ be a $(P_5, C_5, \diam)$-free graph, then it follows from Lemma~\ref{lem:p5.c5.diamond} that $G \in \mathcal{G}_1 \cup \mathcal{G}_2 \cup \mathcal{G}_3$. Since $\mathcal{G}_1$ form a Gallai family, and $\mathcal{G}_2$ is Hamiltonian, the claim holds if $G \in \mathcal{G}_1 \cup \mathcal{G}_2$. Therefore, we may assume that $G \in \mathcal{G}_3$. We let $\mathcal{D}_x$ denote those components of $G\setminus H$ that are complete to only $x$, analogously, $\mathcal{D}_y$ denote those components of $G\setminus H$ that are complete to only $y$. We claim that:

\begin{claim}\label{clm:G3isgalla}
    $\mathcal{G}_3$ is a Gallai family.
\end{claim}

\begin{proof}[Proof of \Cref{clm:G3isgalla}]\renewcommand{\qedsymbol}{}
Let $G$ be any connected graph in $\mathcal{G}_3$. By definition of $\mathcal{G}_3$, let $H = (A,B)$ be the $2P_2$-free bipartite graph contained in $G$, and the vertices $x\in A$, $y\in B$ be such that $xy$ form a dominating edge in $G$. Suppose, without loss of generality, that $d_{H}(x) \geq d_{H}(y)$. We will show that $x$ is a Gallai vertex in $G$. Let $P:=p_1 \dots p_k$ be the longest path in $G$ and assume for the sake of contradiction that $x \notin V(P)$. Note that $V(G\setminus H) \neq \emptyset$, otherwise the claim holds by Theorem~\ref{2p2-free}. Therefore, it follows that either $\mathcal{D}_x$, $\mathcal{D}_y$ are non-empty, or at least one of them is empty;

\begin{itemize}\setlength\itemsep{0.5em}
    \item In the former case, $p_1$ and $p_k$ must be either both in $V(A)$ or one of them must be in $V(\mathcal{D}_y)$, otherwise since $x$ dominates all vertices of $V(B)$ one can include $x$ to $P$ and get a longer path. If $p_1\in V(\mathcal{D}_y)$ and $p_k \in V(A)$, then as $p_{k-1} \in V(B)$ and $p_{k-1}x\in E(G)$, the path $p_1 \dots p_{k-1}xw$ is a longer path for some vertex $w\in \mathcal{D}_x$. The case for $p_1, p_k \in V(A)$ follows by the same argument.

    \item In the latter case, let $\mathcal{D}_x = \emptyset$. If all the vertices in $B$ are contained in $V(P)$, then all the vertices in $B$ are contained in all longest paths of $G$. Otherwise, let $x' \in V(B)$ be such that $x' \notin V(P)$. Then, again as $p_{k-1} \in V(B)$ and $p_{k-1}x\in E(G)$, the path $p_1\dots p_{k-1}xx'$ is a longer path. The case for $\mathcal{D}_y = \emptyset$ is symmetrical.
\end{itemize}

As in both cases we get a contradiction, the vertex $x$ must be common to all longest paths in $G$. This completes the proof of \Cref{clm:G3isgalla}.
\end{proof}

We have shown that the class of $(P_5, C_5, \diam)$-free graphs form a Gallai family. It remains to discuss the case when $G$ contains an induced $C_5$. Then by Lemma~\ref{lem:basicgraphs}, either $G$ is obtained from $G_{1}, \dots, G_{10}$ by stable set expansion of each blue vertex, or $G$ is obtained from $G_{11}$ by stable set expansion of each blue vertex and $P_{3}$-free expansion of the red vertex. For the former case, it is easy to check that $G_{i}$ is $2P_2$-free for $1\leq i \leq 10$ (see Figure~\ref{fig:BS}). Thus, it follows by Theorem~\ref{2p2-free} that $G$ has a Gallai vertex if it is obtained from $G_{1}, \dots, G_{10}$ by a stable set expansion of each blue vertex. For the latter case, we claim that:

\begin{claim}\label{clm:p3-free-expn}
Let $G$ be a class of graphs obtained from $G_{11}$ by a stable set expansion of each blue vertex and $P_{3}$-free expansion of the red vertex. Then $G$ is a Gallai family.
\end{claim}

\begin{proof}[Proof of \Cref{clm:p3-free-expn}]\renewcommand{\qedsymbol}{}
Let $C \subset V(G)$ denote the $P_3$-free expansion of the red vertex. Let $x \in V(G)$ be the vertex complete to all vertices of $C$ (see the graph $G_{11}$ in Figure~\ref{fig:BS}), observe that $V(C) \neq \emptyset$, otherwise $G$ is exactly the graph $G_{2}$ and so we are done by the above discussion. We wish to show that $x$ is a Gallai vertex in $G$. Note that for any vertex $u\in V(G)$, either $u$ is dominated by $x$ or $u$ is dominated by a neighbor of $x$. Let $P:= p_1 \dots p_k$ be a longest path in $G$, and assume for the sake of contradiction that $x \notin V(P)$. Then either of $xp_1$ or $xp_2$ is an edge in $G$. In the former case, $xp_1\dots p_k$ is a path longer than $P$ as $x \notin V(P)$. In the latter case, since $V(C) \neq \emptyset$, there must be a vertex, say $v \in C$, adjacent to $x$ as vertices of $C$ are complete to only $x$. But then $vxp_2 \dots p_k$ is a path longer than $P$, a contradiction. Therefore, $x$ is common to all longest paths in $G$. This proves \Cref{clm:p3-free-expn}.
\end{proof}

Then, \Cref{clm:G3isgalla} combined with \Cref{clm:p3-free-expn} finishes the proof of Theorem~\ref{thm:p5,diam}.
\end{proof}

\begin{observation}\label{alg:P5,diam}
  Let $G$ be a connected $(P_5, \diam)$-free graph. There exists a polynomial-time algorithm which finds a Gallai vertex in $G$.
\end{observation}

\begin{proof}[Proof of Observation~\ref{alg:P5,diam}]{}
There are two possibilities:

\begin{enumerate}
    \item if $G$ has no induced $C_5$, then $G\in \bigcup_{i=1}^{3} \mathcal{G}_i$. If $G \in \mathcal{G}_1 \cup \mathcal{G}_2$, then since $\mathcal{G}_1$ contained in the class of star-like graphs and $\mathcal{G}_2$ is Hamiltonian, it is possible to find all Gallai vertices in $G$ in polynomial time (in fact, for a star-like graph $\mathcal{H}$ with partitions $K, V_1,\dots, V_{k}$ the result of \cite{cerioli2006characterizing} shows that all vertices in $N(v_{i|V_i|}) \cap K$ are Gallai vertices as long as $\mathcal{H}[V_i]$ has the maximum length of the longest path for all $i\in [k]$. Now as $\mathcal{H}[V_i]$ is an interval graph the procedure can be done in polynomial time, see the main result of \cite{ioannidou2009longest}). For when $G \in \mathcal{G}_3$, as we have shown in \Cref{clm:G3isgalla}, the Gallai vertices in $G$ are those of maximum degree among the endpoints $\{x,y\}$ of the dominating edge $e$, which can be found in polynomial time.

    \item If $G$ contains an induced $C_5$, then either $G$ is obtained from $G_1, \dots, G_{10}$ which are $2P_2$-free graphs, or $G$ is obtained from $G_{11}$. In the former case, by the main result of \cite{golanshan2k2}, every vertex of maximum degree is a Gallai vertex. In the latter case, as we have shown in the proof of \Cref{clm:p3-free-expn}, the vertex $x$ which is complete to all vertices of the $P_3$-free expansion set of the red vertex is a Gallai vertex. It is not hard to check that the Gallai vertices can be found in polynomial time in both cases.
\end{enumerate}
This completes the proof.
\end{proof}

\medskip

Finally, we restate and prove \Cref{thm:main}:

\begin{theorem}\label{thm:restateofmain}
Let $H$ be one of the graphs triangle, paw, or diamond. Let $\mathcal{G}$ be the class of graphs restricted to $(P_5, H)$-free graphs, then $\mathcal{G}$ forms a Gallai family. Moreover, there exists a polynomial-time algorithm which finds Gallai vertices on graphs in $\mathcal{G}$.
\end{theorem}

\begin{proof}
The first part of \Cref{thm:restateofmain} follows directly from the proofs of Theorem~\ref{thm:p5,tring}, Corollary~\ref{sta:p5,paw}, and Theorem~\ref{thm:p5,diam}. The second part follows from Observation~\ref{alg:P5,K3} and Observation~\ref{alg:P5,diam}. Note that the algorithm for $(P_5, \paw)$-free graphs is the same as $(P_5, \diam)$-free case. This completes the proof of \Cref{thm:restateofmain}.
\end{proof}

\section{Open problems}

The main open problem related to our work is whether $P_5$-free graphs form a Gallai family. This question was also posed in~\cite{long2023non,golanshan2k2,cerioli2020intersection}. We conjecture the answer to be affirmative. Towards proving it, it would be interesting to consider $(P_5,C_4)$-free graphs and $(P_5,K_4)$-free graphs, which are the last cases of $(P_5,H)$-free graphs with $|H|\leq 4$ that remain open. Note that in the latter case, it is easy to give a constant upper bound on the size of a longest path transversal. Indeed, a result of Bacsó and Tuza~\cite{bacso1990dominating} states that every $P_5$-free graph contains a dominating clique or a dominating $P_3$. By Observation~\ref{2p2-free}, the size of a dominating set upper bounds the size of a longest path transversal. Since in $(P_5,K_4)$-free graphs, a dominating clique has size at most 3, the bound follows.

\bibliographystyle{abbrv}
\bibliography{biblo.bib}

\begin{thebibliography}{10}

\bibitem{arbib2002p5}
C.~Arbib and R.~Mosca.
\newblock On {$(P_5, \text{diamond})$}-free graphs.
\newblock {\em Discrete Mathematics}, 250(1-3):1--22, 2002.

\bibitem{bacso1990dominating}
G.~Bacs{\'o} and Z.~Tuza.
\newblock Dominating cliques in {$P_5$}-free graphs.
\newblock {\em Periodica Mathematica Hungarica}, 21(4):303--308, 1990.

\bibitem{BalisterGLS04}
P.~N. Balister, E.~Gy\H{o}ri, J.~Lehel, and R.~H. Schelp.
\newblock Longest paths in circular arc graphs.
\newblock {\em Combin. Probab. Comput.}, 13(3):311--317, 2004.

\bibitem{CERIOLI2020111717}
M.~R. Cerioli, C.~G. Fernandes, R.~Gómez, J.~Gutiérrez, and P.~T. Lima.
\newblock Transversals of longest paths.
\newblock {\em Discrete Mathematics}, 343(3):111717, 2020.

\bibitem{cerioli2020intersection}
M.~R. Cerioli and P.~T. Lima.
\newblock Intersection of longest paths in graph classes.
\newblock {\em Discrete Applied Mathematics}, 281:96--105, 2020.

\bibitem{cerioli2006characterizing}
M.~R. Cerioli and J.~L. Szwarcfiter.
\newblock Characterizing intersection graphs af substars of a star.
\newblock {\em Ars Combinatoria}, 79:21--31, 2006.

\bibitem{Chen15}
F.~Chen.
\newblock Nonempty intersection of longest paths in a graph with small matching
  number.
\newblock {\em Czechoslovak Mathematical Journal}, 65:545--553, 2015.

\bibitem{ChenEFHSYY17}
G.~Chen, J.~Ehrenm\"uller, C.~G. Fernandes, C.~G. Heise, S.~Shan, P.~Yang, and
  A.~N. Yates.
\newblock Nonempty intersection of longest paths in series-parallel graphs.
\newblock {\em Discrete Mathematics}, 340(3):287--304, 2017.

\bibitem{choudum2010first}
S.~A. Choudum and T.~Karthick.
\newblock First-fit coloring of {$\{P_5, K_{4}-e\}$}-free graphs.
\newblock {\em Discrete Applied Mathematics}, 158(6):620--626, 2010.

\bibitem{Gallai68}
P.~Erd\H{o}s and G.~Katona, editors.
\newblock {\em Theory of Graphs}.
\newblock Proceedings of the Colloquium held at Tihany, Hungary, September
  1966. Academic Press, New York, 1968.
\newblock Problem 4 (T. Gallai), p. 362.

\bibitem{gao2021nonempty}
Y.~Gao and S.~Shan.
\newblock Nonempty intersection of longest paths in graphs without forbidden
  pairs.
\newblock {\em Discrete Applied Mathematics}, 304:76--83, 2021.

\bibitem{golanshan2k2}
G.~Golan and S.~Shan.
\newblock Non-empty intersection of longest paths in {$2K_2$}-free graphs.
\newblock {\em The Electronic Journal of Combinatorics}, 25:P2.37, 2018.

\bibitem{ioannidou2009longest}
K.~Ioannidou, G.~B. Mertzios, and S.~D. Nikolopoulos.
\newblock The longest path problem is polynomial on interval graphs.
\newblock In {\em Mathematical Foundations of Computer Science 2009: 34th
  International Symposium, MFCS 2009, Novy Smokovec, High Tatras, Slovakia,
  August 24-28, 2009. Proceedings 34}, pages 403--414. Springer, 2009.

\bibitem{JobsonKLW16}
A.~S. Jobson, A.~E. K\'ezdy, J.~Lehel, and S.~C. White.
\newblock Detour trees.
\newblock {\em Discrete Applied Mathematics}, 206:73--80, 2016.

\bibitem{Joos15}
F.~Joos.
\newblock A note on longest paths in circular arc graphs.
\newblock {\em Discussiones Mathematicae Graph Theory}, 35(3):419--426, 2015.

\bibitem{sublinearLPT}
J.~A. Long, K.~G. Milans, and A.~Munaro.
\newblock Sublinear longest path transversals.
\newblock {\em SIAM Journal on Discrete Mathematics}, 35(3):1673--1677, 2021.

\bibitem{long2023non}
J.~A. Long~Jr, K.~G. Milans, and A.~Munaro.
\newblock Non-empty intersection of longest paths in {$H $}-free graphs.
\newblock {\em The Electronic Journal of Combinatorics}, page P1.32, 2023.

\bibitem{olariu1988paw}
S.~Olariu.
\newblock Paw-fee graphs.
\newblock {\em Information Processing Letters}, 28(1):53--54, 1988.

\bibitem{rautenbach2014}
D.~Rautenbach and J.-S. Sereni.
\newblock Transversals of longest paths and cycles.
\newblock {\em SIAM Journal on Discrete Mathematics}, 28(1):335--341, 2014.

\bibitem{sumner1981subtrees}
D.~P. Sumner.
\newblock Subtrees of a graph and chromatic number.
\newblock {\em The Theory and Applications of Graphs,(G. Chartrand, ed.), John
  Wiley \& Sons, New York}, 557:576, 1981.

\bibitem{walther1969nichtexistenz}
H.~Walther.
\newblock {\"U}ber die nichtexistenz eines knotenpunktes, durch den alle
  l{\"a}ngsten wege eines graphen gehen.
\newblock {\em Journal of Combinatorial Theory}, 6(1):1--6, 1969.

\bibitem{zamfirescu1976longest}
T.~Zamfirescu.
\newblock On longest paths and circuits in graphs.
\newblock {\em Mathematica Scandinavica}, 38(2):211--239, 1976.

\end{thebibliography}
\end{document}